\documentclass[12pt]{amsart}
\usepackage{amsmath,amssymb,amscd,array, mathrsfs }

\usepackage{amsmath,amscd,amssymb,amsthm,array}
\usepackage{color}

\usepackage{amsmath,amssymb,mathrsfs,amsthm, tikz-cd,mathrsfs}
\usepackage{comment}
\def\url#1{\expandafter\s

\tring\csname #1\endcsname}

\def\mmat #1,#2,#3,#4,{\text{\small\arraycolsep=3pt $
\begin{pmatrix}#1&#2\\#3&#4\end{pmatrix}$}}

\usepackage{hyperref}
\usepackage{capt-of}

\usepackage{comments}
\newComments\SBe{Said}{blue}
\newComments\SBo{Sofiane}{blue}
\newComments\AM{Nacer}{blue}
\newComments\DL{DL}{red}
\newComments\QEh{QEh}{blue}

\def\mmat #1,#2,#3,#4,{\text{\small\arraycolsep=3pt $
\begin{pmatrix}#1&#2\\#3&#4\end{pmatrix}$}}

\usepackage{lscape}
\usepackage{tikz-cd}
\usepackage{enumerate}

\usepackage{DLdef1}
\usepackage{multicol}

\def\mmat #1,#2,#3,#4,{\text{\small\arraycolsep=3pt $
\begin{pmatrix}#1&#2\\#3&#4\end{pmatrix}$}}

\renewcommand {\ssbegin}[2][*]
 {\refstepcounter{subsection}%
\if#1*
\addcontentsline{toc}{subsection}{\thesubsection.\hskip 1pc #2}%
\else
\addcontentsline{toc}{subsection}{\thesubsection.\hskip 1pc #2. #1}%
\fi
 \def \secno {\gdef \secno {}{\ssecfont
\thesubsection.\hskip 2ex}%
 }%
 \begin{#2}}

\renewcommand {\sssbegin}[2][*]
 {\refstepcounter{subsubsection}
\if#1*
\addcontentsline{toc}{subsubsection}{\thesubsubsection.\hskip 1pc #2}%
\else
\addcontentsline{toc}{subsubsection}{\thesubsubsection.\hskip 1pc #2. #1}
\fi
 \def \secno {\gdef \secno {}{\ssecfont \thesubsubsection.\hskip 2ex}%
 }%
 \begin{#2}}

\renewcommand {\parbegin}[2][*]
 {\refstepcounter{paragraph}
\if#1*
\addcontentsline{toc}{paragraph}{\theparagraph.\hskip 1pc #2}%
\else
\addcontentsline{toc}{paragraph}{\theparagraph.\hskip 1pc #2. #1}
\fi
 \def \secno {\gdef \secno {}{\ssecfont \theparagraph.\hskip 2ex}%
 }%
 \begin{#2}}

\newcommand {\ce}{{\text{CE}}}

\newcommand {\fj}{{\mathfrak{j}}}

\rmnameii{etr}{etr}
\rmnameii{evv}{ev}

\DeclareMathOperator{\K}{\mathbb{K}}

\DeclareMathOperator{\sh}{\text{Sh}}

\newcommand {\w}{\omega}

\newcommand{\Z}{\mathbb{Z}}

\begin{document}

\title[Double extensions of quasi-Frobenius Lie superalgebras]{Double extensions of quasi-Frobenius Lie superalgebras with degenerate center}

\author{Sofiane Bouarroudj}

\address{Division of Science and Mathematics, New York University Abu Dhabi, P.O. Box 129188, Abu Dhabi, United Arab Emirates.}
\email{sofiane.bouarroudj@nyu.edu}

\author{Quentin Ehret}
\address {Division of Science and Mathematics, New York University Abu Dhabi, P.O. Box 129188, Abu Dhabi, United Arab Emirates.}
\email{qe209@nyu.edu}

\thanks{SB and QE were supported by the grant NYUAD-065.}






\begin{abstract}
We develop the process of symplectic double extensions for Lie superalgebras with degenerate center. The construction is a superization of a recent work by Fischer, and generalize our previous work. We provide a standard model for such double extensions, where the symplectic form is either orthosymplectic or  periplectic. Additionally, we show  that every double extension is naturally equivalent to either of these two standard types of extensions. Several examples in low dimensions are given to illustrate the concept.   
\end{abstract}


\maketitle

\thispagestyle{empty}
\setcounter{tocdepth}{2}
\tableofcontents

\newpage
\section{Introduction} \label{SecDef}
\subsection{Double extensions of quadratic Lie (super)algebras} Let $\fg$ be a Lie algebra with a non-degenerate {\it symmetric} bilinear form (called {\it quadratic} in the literature). The notion of {\it  double extension} of quadratic Lie algebras have been fully studied by Medina and Revoy in \cite{MR1}. They showed that any irreducible non-simple quadratic Lie algebra can be obtained as a double extension. In particular, in the case where $\fg$ admits a non-trivial center, then $\fg$ can be obtained as a double extension of a smaller algebra by a $1$-dimensional space. Similar results were also obtained in \cite{FS} and a classification in low dimension is also provided. 

A superization of the results of \cite{MR1} was first explored in \cite{BeBe} when the odd part is $2$-dimensional. Since useful tools such as the Lie Theorem and the Levi decomposition do not hold in general for Lie superalgebras, new ideas are needed. Those results have been improved in \cite{BBB}, and in subsequent papers, see \cite{B} and references therein.

 Over a field of characteristic $p\not = 0 $, the process of double extensions is still valid, except in characteristic two where new techniques are needed due to the presence of the squaring, see \cite{BB1,BB2, BLS}.

 The study of double extensions of \textit{restricted} Lie superalgebras started in \cite{BBH}, where quadratic Lie superalgebras were called {\it NIS} (an abbreviation for `non-degenerate invariant symmetric'). Under certain conditions, it is possible to endow a double extension with a restricted structure. There are also several conditions under which a quadratic restricted superalgebra can be obtained as a double extension. The conditions arise from the restricted cohomology, see \cite{BBH, BEM}. 

\subsection{Double extensions of quasi-Frobenius Lie algebras} Roughly speaking, a Lie algebra is called quasi-Frobenius (resp. Frobenius) if it is equipped with a non-degenerate closed (resp. exact) antisymmetric bilinear form, called \textit{symplectic}. Frobenius Lie algebras were introduced by Ooms in \cite{O} to answer a question raised by Jacobson.  Double extensions of quasi-Frobenius Lie algebras (by a symplectic plane) were introduced by Medina and Revoy in \cite{MR2} in order to classify nilpotent symplectic groups by their algebras. Indeed, there is a one-to-one correspondence between quasi-Frobenius Lie algebras and simply connected Lie groups admitting a left-invariant symplectic form. Medina and Revoy showed that any nilpotent quasi-Frobenius Lie algebra can be obtained by a sequence of double extensions from the trivial algebra.
In addition, Baues and Cort\'es conducted an extensive study in \cite{BC}, where they introduced the concepts of \textit{symplectic reduction} and \textit{symplectic oxidation}, which reverse each other. The symplectic oxidation process is closely related to double extensions. In particular, they applied their constructions to prove that any symplectic Lie group admits a sequence of symplectic reductions that leads to a unique irreducible symplectic Lie group. 

\subsection{On the existence of  isotropic ideals} Let $(\fg,\w)$ be a symplectic Lie algebra and let $\fj\subseteq\fg$ be an {\it isotropic} ideal, see Section \ref{2.3}. A symplectic reduction of $(\fg,\w)$ is a symplectic structure on the quotient  Lie algebra $\fj^{\perp}/\fj$. It is crucial to determine an appropriate isotropic ideal $\fj$, not just for reduction but also for double extension. In \cite{BC}, 
 the authors provide a procedure to find such an ideal in some particular cases (see \cite[Thms. 8.3 and 8.5]{BC}). They also give examples of 10-dimensional algebras without Lagrangian ideals.  In \cite[Lemme 1.1]{DM2}, Dardi\'e and Medina proved that any symplectic Lie algebra admits a non trivial abelian ideal $\fj$, which is either isotropic or symplectic. Moreover, they showed that in the case where a symplectic algebra $\fg$ does not admit any isotropic ideal, then $\fg$ is the semi-direct product $[\fg,\fg]\rtimes[\fg,\fg]^{\perp}$ (\cite[Th\'eor\`eme 1.3]{DM2}).
 
In the case where $\fg$ admits a degenerate center (for instance, nilpotent non-abelian Lie algebras have a degenerate center), Fisher has developed a method for selecting the ideal $\fj$ canonically, see \cite{Fi}. In particular, his method makes it possible to consider double extensions by ideals of arbitrary dimensions. Moreover, he constructed a standard model for these extensions and showed that any double extension is equivalent to this standard model. This result was used to prove that a quasi-Frobenius Lie algebra with a degenerate center is obtained by means of such an extension. Further, he performed a classification in dimension 6 for nilpotent algebras.
This double extension process cannot be performed in general in the case where the center of the algebra is trivial. Bajo and Benayadi overcame this difficulty in the case where the Lie algebra is abelian para-K\"ahler in \cite{BaBe} and were able to describe double extensions even if the center is trivial. The case where $\fg$ is solvable but not necessarily abelian para-K\"ahler has not yet been solved.

Complex symplectic oxidation (double extension) was introduced by Bazzoni et al. in \cite{BFLM}, where they obtained all nilpotent complex symplectic Lie algebras in dimension 8.

\subsection{Double extensions of quasi-Frobenius Lie superalgebras}
The results of \cite{MR2} have been superized in \cite{BM}, where the authors divided the study into four cases, depending on the parity of the symplectic form and of the derivation used to build the double extension. In particular, they showed that every quasi-Frobenius Lie superalgebra can be realized as double extension by a symplectic plane under some conditions and provided a classification of quasi-Frobenius Lie superalgebras in dimension $4$. The study has been carried out later to the case of restricted Lie superalgebras in \cite{BEM}. The main novelty of this work is that restricted cohomology captures the obstructions to symplectic double extensions.

As this paper follows Fisher's ideas (\cite{Fi}), we will look at the challenges involved in the super settings. We prove that choosing canonically the ideal $\fj$ leads a superization of Fisher's work. There are two cases to consider: either the symplectic form is even or odd. It is necessary that the symplectic form is homogeneous; otherwise, the ideal $\fj^\perp$ is not homogeneous. 

The paper is organized as follows. Basic concepts about Lie superalgebras, quasi-Frobenius structures and cohomological tools are introduced in Section \ref{background}. The core of the paper is contained is Section \ref{sectiondoubleext}, where we introduce symplectic double extensions of a quasi-Frobenius Lie superalgebra by an abelian Lie superalgebra (Definition \ref{doublextdef}) as well as equivalence of such extensions. We then provide a standard model in the orthosymplectic case (Section \ref{subsectionortho}) and the periplectic case (Section \ref{subsectionperi}). The main result of this paper is Theorem \ref{equivstandard} which states that any double extension is equivalent to the appropriate standard model.  We provide in Section \ref{sectionexamples} a list of explicit examples of double extensions based on the classification done in \cite{Ba}.\\

\noindent\textbf{Conventions and notations.} The group of integers modulo $2$ is denoted by $\Z_{2}$. Let $V$ be a $\Z_{2}$-graded space. The degree of an homogeneous element $v\in V_{\bar{i}}$ is denoted by $|v|:=\bar{i}$. The element $v$ is called \textit{even} if $v\in V_\ev$ and \textit{odd} if $v\in V_\od$. A linear map $\varphi:V\rightarrow W$ between $\Z_{2}$-graded space is called \textit{even} if $\varphi(V_{\bar{i}})\subset W_{\bar{i}}$ and \textit{odd} if $\varphi(V_{\bar{i}})\subset W_{\overline{i+1}}$.\\~\\
Throughout the paper, $\K$ will be a field of characteristic $0$. 

\section{Background}\label{background} 
A comprehensive introduction to Lie superalgebras can be found in \cite{L, S}.

\subsection{Lie superalgebras}

A \textit{Lie superalgebra} is a $\Z_2$-graded vector space $\fg=\fg_\ev\oplus \fg_\od$ equipped with a bilinear map $[\cdot,\cdot]:\fg\times\fg\rightarrow\fg$ such that, for all homogeneous elements $x,y,z\in\fg$, we have

\begin{enumerate}[($i$)]
    \item $[\fg_{\bar{i}},\fg_{\bar{j}}]\subseteq \fg_{\overline{i+j}}$;
    \item $[x,y]=-(-1)^{|x||y|}[y,x]$;
    \item $(-1)^{|x||z|}[[x,y],z]+(-1)^{|x||y|}[[y,z],x]+(-1)^{|y||z|}[[z,x],y]=0$ (Jacobi identity).
\end{enumerate}
From the definition $\fg_\ev$ is an ordinary Lie algebra, and $\fg_\od$ is a $\fg_\ev$-module. 

A linear endomorphism  $d:\fg\rightarrow\fg$ is called  a \textit{derivation} if it satisfies
$$d([x,y])=[d(x),y]+(-1)^{|d||x|}[x,d(y)],~\forall x,y\in \fg.$$   
In fact, this condition is cohomological. It means that $d$ is a 1-cocyle on $\fg$ with values in the adjoint representation. 

An \textit{ideal} of $\fg$ is a subspace $\fj\subseteq\fg$ satisfying $[\fj,\fg]\subseteq\fj$. it is called \textit{homogeneous} if, in addition,  $\fj=(\fj\cap\fg_\ev)\oplus (\fj\cap\fg_\od).$

The  \textit{center} of $\fg$ is the homogeneous ideal $\fz:=\{z\in\fg,~[z,x]=0~\forall x\in\fg\}$.

\subsection{Change of parity functor.} Let $V=V_\ev\oplus V_\od$ a $\Z_2$-graded vector space. We denote by $\Pi$ the \textit{change of parity functor} $\Pi: V\mapsto \Pi (V)$, where $\Pi V$ is another copy of $V$ such that $ \Pi(V)_\ev:=V_\od;~~\Pi(V)_\od:=V_\ev$ (see \cite{L}).  Elements of $\Pi(V)$ shall be denoted by $\pi(v),~\forall v\in V$. We also shall identify $V$ and $\Pi(\Pi(V))$ in a natural way.

\subsection{Quasi-Frobenius Lie superalgebras}\label{2.3}    
Let $\fg$ be a Lie superalgebra and let $\w$ be a bilinear form on $\fg$. Recall that the form $\w$ is called \textit{non-degenerate} if
$$\ker(\w):=\{x\in\fg,~\omega(x,y)=0~\forall y\in\fg\}=\{0\}.$$
Following \cite{BaBe, BM}, a Lie superalgebra $\fg$ is called {\it quasi-Frobenius} if it is equipped with a non-degenerate 2-cocycle $\omega\in Z^2_ {\text{CE}}(\fg; \Kee)$ (see Section \ref{cohomconstr}), that is
\begin{equation}\label{coccond}
(-1)^{|a||c|} \omega(a, [b , c]) +(-1)^{|c||b|} \omega(c, [a , b] )+(-1)^{|b||a|} \omega(b, [c , a] )=0.
\end{equation}
 We denote such an algebra by $(\fg, [\cdot,\cdot], \omega)$. The bracket will be omitted whenever it is clear from the context. 
 
 In \cite{BM}, we described filliform Lie superalgebras as quasi-Frobenius and list quasi-Frobenius Lie superalgebras in small dimensions.
 
 A bilinear form $\w$ satisfying \eqref{coccond} is called \textit{closed}. In the case where $\omega \in B^2(\fg,\mathbb K)$, the Lie superalgebra $\fg$ is called  {\it Frobenius} while $\w$ is called \textit{exact}.
 
A quasi-Frobenius Lie superalgebra $(\fg, \omega)$ is called 

\begin{enumerate}[$(i)$]
    \item \textit{orthosymplectic} if the form $\w$ is even, that is, $\w(\fg_\ev,\fg_\od)=0$;
    \item \textit{periplectic} if the form $\w$ is odd, that is, $\w(\fg_\ev,\fg_\ev)=\w(\fg_\od,\fg_\od)=0$.
\end{enumerate}
Let $\fj\subseteq\fg$ be a subspace of a quasi-Frobenius Lie superalgebra $(\fg,\w)$. The \textit{orthogonal} of $\fj$ (with respect to the form $\w$) is defined by
$$\fj^{\perp}:=\{x\in\fg,~\w(x,y)=0~\forall y\in\fj \}.  $$
An ideal $\fj\subseteq\fg$ is called 

\begin{enumerate}[$(i)$]

\item \textit{isotropic} if $\fj\subset\fj^{\perp}$;
\item \textit{lagrangian} if $\fj$ is maximal and if $\fj=\fj^{\perp}$; 
\item \textit{degenerate} if $\fj\cap\fj^{\perp}\neq \{0\}$.
\end{enumerate}

A morphism of quasi-Frobenius Lie superalgebras $\varphi:(\fg,\w_\fg)\rightarrow(\fh,\w_\fh)$ is a morphism of Lie superalgebras satisfying $\varphi^*\w_\fg=\w_\fh$, where $\varphi^*\w_\fg(x,y):=\w_{\fh}(\varphi(x),\varphi(y))$, for all $x,y\in\fg$.

Quasi-Frobenius Lie (super)algebras are also called symplectic Lie (super)algebras, see \cite{BC, D, DM1, DM2, Fi, MR2, RS}. 

\sssbegin{Lemma}[Superization of \text{\cite[Lemma 2.1]{BC}}]
    Let $(\fg,\omega)$ be a quasi-Frobenius Lie superalgebra and $\fj$ be an isotropic homogeneous ideal of $\fg$. Then $\fj$ is abelian.
\end{Lemma}

\begin{proof}
Let $x,y\in \fj$. Since $\fj\in\fj^{\perp}$, we get $\omega([x,y],z)=0,~\forall z\in \fg$. It follows that $[x,y]=0$.
\end{proof}

\sssbegin{Proposition}\label{quotient}
    Let $(\fg,\omega)$ be a quasi-Frobenius Lie superalgebra such that $\w_{\fg}$ is homogeneous and let $\fj$ be an isotropic homogeneous ideal of $\fg$. Then $\overline{\fg}:=\fj^{\perp}/\fj=\fj^{\perp}_\ev/\fj_\ev\oplus \fj^{\perp}_\od/\fj_\od$ is a quasi-Frobenius Lie superalgebra. Moreover, $|\w_{\fg}|=|\w_{\overline{\fg}}|$.
\end{Proposition}

\begin{proof}
 Observe that $\fj^{\perp}$ is homogeneous. Indeed, let $x=x_\ev+x_\od\in \fj^{\perp}$ and $y=y_\ev+y_\od\in \fj$. From $\w(x,y_\ev)=0$, we obtain that $\w(x_\ev,y_\ev)=\w(x_\od,y_\ev)=0$. From $\w(x,y_\od)=0$, we obtain that $\w(x_\ev,y_\od)=\w(x_\od,y_\od)=0$. We deduce that $x_{\ev}\in \fj^{\perp}$ and that $x_{\od}\in \fj^{\perp}$. Therefore, $\fj^{\perp}$ is homogeneous. Since $\fj$ is homogeneous, the quotient $\overline{\fg}=\fj^{\perp}/\fj$ is a Lie superalgebra.

We define a symplectic form on $\overline{\fg}$ by $\w_{\overline{\fg}}(\overline{x},\overline{y}):=\w_{\fg}(x,y)$. The form $\w_{\overline{\fg}}$ is well defined, non-degenerate and closed. Therefore, the Lie superalgebra $(\overline{\fg},\w_{\overline{\fg}})$ is quasi-Frobenius.
\end{proof}

\subsection{Cohomological constructions}\label{cohomconstr}
Let $\fg$ be a Lie superalgebra and $M$ be an $\fg$-module. Let us recall the space of cochains $C^n_{\text{CE}}(\fg; M)$ in the Chevalley-Eilenberg cohomology. For $n=0$, we put $C^0_{\text{CE}}(\fg; M):=M$. For $n>0$, the space of cochains $C^n_{\text{CE}}(\fg; M)$ is the space of $n$-linear super anti-symmetric maps. The differential is given as follows:

\[
\begin{array}{l}
 d_{\mathrm{CE}}^0(m)(x)=(-1)^{|m||x|}x\cdot m\text{~~for any $m\in M$ and $x\in\fg$};\\[2mm]
d_{\mathrm{CE}}^n(\varphi)(x_1,\cdots, x_n) \\
= \mathop{\sum}\limits_{i<j}(-1)^{(|x_j|(|x_{i+1}|+\dots+|x_{j-1}|)+j}\varphi(x_1, \cdots, x_{i-1}, [x_i,x_j], x_{i+1} \cdots, \widetilde  x_j, \cdots, x_n)\\
+\mathop{\sum}\limits_{j}(-1)^{|x_j|(|f|+|x_{1}|+\cdots+|x_{j-1}|)+j}x_j \cdot \varphi (x_1, \cdots, \widetilde  x_j, \cdots, x_n)\\
\text{~~for any $\varphi \in C^{n-1}_{\text{CE}}(\fg; M)$ with $n>0$, and $x_1,\dots, x_n\in\fg$}. \\
\end{array}
\]
Let us denote the spaces $Z^n_{\mathrm{CE}}(L;M):=\Ker\left(d^n_{\mathrm{CE}}\right)$ and $B^n_{\mathrm{CE}}(L;M):=\text{Im}\left(d^{n-1}_{\mathrm{CE}}\right)$, for all $n\geq1$. The quotient space $\mathrm{H}^n_{\mathrm{CE}}(L;M):=Z^n_{\mathrm{CE}}(L;M)/B^n_{\mathrm{CE}}(L;M)$ is well defined for all $n\geq 0$ and is called the $n^{th}$ \textit{Chevalley-Eilenberg cohomology space of} $\fg$ with coefficients in the $\fg$-module $M$.

Let $n,k$ be natural integers. We denote by $\sh(n,k)$ the set of all permutations $\tau$ of $\left\lbrace 1,\cdots,n,\cdots,n+k \right\rbrace$ such that $\tau(1)<\tau(2)<\cdots<\tau(n) \text{ and } \tau (n+1)<\cdots<\tau(n+k).$
Let $V=V_\ev\oplus V_\od$ be a $\Z_2$-graded space and let $X=\{x_1,\cdots,x_n\}$ be a subset of homogeneous elements of $V$. Let $\sigma$ be a permutation of $\left\lbrace 1,\cdots,n,\cdots,n+k\right\rbrace$. Following \cite{BP}, we denote by
$$K(\sigma,X):=\text{Card}\bigl\{(i,j),~i<j,~x_{\sigma(i)}\in V_\od,~x_{\sigma(j)}\in V_\od,~\sigma(i)>\sigma(j)\bigl\}.$$ Then, we denote by $\theta(\sigma,X):=\text{sgn}(\sigma)(-1)^{K(\sigma,X)}$, where $\text{sgn}(\sigma)$ is the signature of the permutation $\sigma$.

Let $\fg$ be a Lie superalgebra and let $U,V,W$ be $\fg$-modules. Recall that a bilinear map $m:U\times V\rightarrow W$ is called $\fg$-equivariant if
\[
x\cdot m(u,v)= (-1)^{|x||m|}m(x \cdot u,v)+(-1)^{|x|(|u|+|m|)}m(u,x\cdot v),~\forall x\in\fg,~\forall (u,v)\in U\times V \text{ homogeneous.}
\]
The following definition is crucial to us. 
\sssbegin{Definition}
    Let $\fg$ be a Lie superalgebra, let $U,V,W$ be $\fg$-modules and let\\ $m:U\times V\rightarrow W$ be a bilinear $\fg$-equivariant map. We define a new product  $C^n_{\ce}(\fg,U)\times C^k_{\ce}(\fg,V)\rightarrow C^{n+k}_{\ce}(\fg,W)$ by
    \begin{equation}
        m(\alpha\wedge\beta)(x_1,\cdots,x_{n+k}):=\sum_{\sigma\in\sh(n,k)}\theta(\sigma,X)m\bigl(\alpha(x_{\sigma(1)},\cdots,x_{n}),\beta(x_{\sigma(n+1)},\cdots,x_{n+k})\bigl).
    \end{equation}
\end{Definition}

Let $\fa$ be a Lie superalgebra and let $\fh$ be a trivial $\fa$-module. We define an evaluation map
\begin{equation}\label{evaluation}
    \Ev:C^1_{\ce}(\fa,\fh)\times \fa\rightarrow\fh,~(f,a)\longmapsto f(a).
\end{equation}
Let $\xi:\fa\rightarrow\Hom(\fh)$ be a linear map. We define a map $d_{\xi}:C^n_{\ce}(\fa,\fh)\rightarrow C^{n+1}_{\ce}(\fa,\fh)$ by
\begin{align}
    d_{\xi}\varphi(a_1,\cdots,a_{n+1})=&\sum_{i<j}(-1)^{(|a_j|(|a_{i+1}|+\dots+|a_{j-1}|)+j}\varphi(a_1, \dots, a_{i-1}, [a_i,a_j], a_{i+1} \dots, \widetilde  a_j, \dots, a_{n+1})\nonumber\\
    &+\sum_{j}(-1)^{|a_j|(|\varphi|+|a_{1}|+\dots+|a_{j-1}|)+j}\xi(x_j) \bigl( \varphi (a_1, \dots, \widetilde  a_j, \dots, a_{n+1})\bigl).
\end{align}

\noindent\textbf{Remark.} The map $d_{\xi}$ satisfies $d_{\xi}^{n+1}\circ d_{\xi}^n=0$ if and only if $\xi$ is a Lie superalgebras map.

These objects were used for example in \cite{N} to describe non-abelian extensions of topological Lie algebras. We will need them in our constructions, see e.g. Thms. \ref{orthostandard} and \ref{peristandard}.

\section{Symplectic double extensions}\label{sectiondoubleext} 
The following definition is a straightforward superization of \cite{Fi}. 

\ssbegin{Definition}[Symplectic double extension]\label{doublextdef}
Let $(\fa,[\cdot,\cdot]_\fa,\w_{\fa})$ be a quasi-Frobenius Lie superalgebra and let $\fl$ be an abelian Lie superalgebra. Following \cite{Fi}, a symplectic double extension of $\fa$ by $\fl$ is given by a quadruple $(\fg,\fj,i,p)$, where
\begin{enumerate}[($i$)]
    \item $\fg=(\fg,[\cdot,\cdot]_{\fg},\w_\fg)$ is a quasi-Frobenius Lie superalgebra such that $\w_{\fg}$ is homogeneous;
    \item $\fj$ is a central isotropic homogeneous ideal of $\fg$ (with respect to $\w_\fg$);
    \item There is a short exact sequence of Lie superalgebras
        \begin{equation}
            0\longrightarrow\fa\overset{i}{\longrightarrow}\fg/\fj\overset{p}{\longrightarrow}\fl\longrightarrow 0,
        \end{equation} where $i:\fa\rightarrow \fj^{\perp}/\fj$ is an isomorphism of quasi-Frobenius Lie superalgebras. Following \cite{Fi}, we call such an extension \textit{balanced} if $\fj=\fz\cap\fz^{\perp}$, where $\fz$ is the center of $\fg$.
\end{enumerate}
\end{Definition}

Actually, any central ideal $\fj$ is homogeneous. Therefore, the quotient $\fj^{\perp}/\fj$ is a quasi-Frobenius Lie superalgebra by Prop. \ref{quotient}.

\sssbegin{Theorem}[Superization of \text{\cite[Theorem 1]{Fi}}]
    Let $(\fg,\omega_\fg)$ be a quasi-Frobenius Lie superalgebra with degenerate center $\fz$. Suppose that $\omega_\fg$ is homogeneous. Then $\fz^{\perp}$ is homogeneous and $(\fg,\omega_\fg)$ has the structure of a balanced quadratic extension in a canonical way.
\end{Theorem}

\begin{proof}
    The center $\fz$ is always homogeneous. Since $\w_{\fg}$ is closed, $\fz^{\perp}$ is an homogeneous ideal. Therefore, $\fz\cap\fz^{\perp}$ is homogeneous. Choose $\fj=\fz\cap\fz^{\perp}$. The proof follows as in \cite[Theorem 1]{Fi}.
\end{proof}

\sssbegin{Definition}[Equivalence of symplectic double extensions]\label{eqextt}
Let $(\fg_1,\fj_1,i_1,p_1)$ and $(\fg_2,\fj_2,i_2,p_2)$ be two symplectic double extensions of a quasi-Frobenius Lie superalgebra $(\fa,\w_{\fa})$ by an abelian Lie superalgebra $\fl$. An equivalence of symplectic double extensions is a symplectic isomorphism $\phi:\fg_1\rightarrow\fg_2$ satisfying $\phi(\fj_1)=\fj_2$ and such that the following diagram commutes:
\begin{equation}\label{eqext}
\begin{tikzcd}
            &                                                & \fg_1/\fj_1\arrow[rd, "p_1"] \arrow[dd, "\overline{\phi}"] &             &   \\
0 \arrow[r] & \fa \arrow[ru, "i_1"] \arrow[rd, "i_2"'] &                                              & \fl \arrow[r] & 0, \\
            &                                                & \fg_2/\fj_2 \arrow[ru, "p_2"']                     &             &  
\end{tikzcd}
\end{equation} where $\overline{\phi}$ is the induced map between the quotient spaces.
\end{Definition}

Let $a\in\fa$. From the definition of the equivalence, we have
$$\phi\circ i_1(a)+\fj_2=\phi(i_1(a)+\fj_1)=\overline{\phi}\circ i_1(a)=i_2(a)=a+\fj_2,~\forall a\in\fa$$ Therefore, $\phi\circ i_1(a)-a\in\fj_2,~\forall a\in\fa$.

\subsection{A standard model for double extensions -- the orthosymplectic case }\label{subsectionortho}

 Let $(\fa,[\cdot,\cdot]_\fa,\omega_\fa)$ be an orthosymplectic quasi-Frobenius Lie superalgebra and let $\fl$ be an abelian Lie superalgebra. The dual of $\fl$ is denoted by $\fl^*$. We define a bracket $[\cdot,\cdot]_\fd$ on $\fd:=\fl^*\oplus\fa\oplus\fl$ as follows: 
\begin{equation}\label{extbracket}
\begin{array}{lclrcll}
[\fl^*,\fd ] & : = & 0; &   & [a,b]_\fd&:=&\beta(a,b)+[a,b]_\fa; \\[2mm] 
 [L,a]_\fd & := & \gamma(L)(a)+\xi(L)(a);& & [L_1,L_2]_\fd&:=&\epsilon(L_1,L_2)+\alpha(L_1,L_2),
 \end{array}
\end{equation}
for all $a,b\in\fa$ and for all $L,L_1,L_2\in\fl$, where
\begin{align}
    \label{beta} 
 \beta:&~\fa\wedge\fa\rightarrow\fl^*,(a,b)\mapsto\beta(a,b)(L):=-(-1)^{|L|(|a|+|b|)}\omega_\fa\bigl(\xi(L)(a),b\bigl)\\\nonumber&\hspace{4cm}-(-1)^{|b||L|}\omega_\fa\bigl(a,\xi(L)(b)\bigl),~\forall L\in \fl;\\
    \xi:&~\fl\rightarrow\End(\fa) \text{ is linear;}\\
    \gamma:&~\fl\rightarrow\Hom(\fa,\fl^*) \text{ is linear;}\\
    \epsilon:&~\fl\times\fl\rightarrow\fl^* \text{ is linear;}\\
    \alpha:&~\fl\times\fl\rightarrow\fa \text{ is bilinear and defined by}\\\label{alpha}
    \omega_{\fa}\bigl(a,\alpha(L_1,L_2)\bigl)&:=(-1)^{|L_2|(|L_1|+|a|)}\gamma(L_2)(a)(L_1)\\\nonumber&~~-(-1)^{|a||L_1|}\gamma(L_1)(a)(L_2),~\forall a\in \fa,~\forall L_1,L_2\in \fl.  
\end{align}

Moreover, we define the following induced maps:
\begin{align}
    \forall a\in\fa,~\widetilde{\beta}(a)&:=\beta(a,\cdot):\fa\rightarrow\fl^*;\\
    \forall a\in\fa,~\forall L\in\fl,~\widetilde{\gamma}(a)(L)&:=(-1)^{|a||L|}\gamma(L)(a)\in \fl^*;\\
        \forall a\in\fa,~\forall L\in\fl,~\widetilde{\xi}(a)(L)&:=(-1)^{|a||L|}\xi(L)(a)\in\fa;\\
        \forall L\in\fl,~\forall f\in C^1_{\ce}(\fa,\fl^*),~\widehat{\xi}(L)(f)&:=(-1)^{|f||L|}f\circ\xi(L):\fa\rightarrow\fl^*.    
\end{align}

\noindent\textbf{Remark.} In the case where $\fl$ is one-dimensional, the map $\beta$ defined in Eq. \eqref{beta} agrees with the map $\phi$ as defined in \cite[Theorem 4.1.1]{BM}. 

\sssbegin{Theorem}[Superization of \text{\cite[Theorem 2]{Fi}}]\label{orthostandard}
Let $(\fa,[\cdot,\cdot]_\fa,\omega_\fa)$ be an orthosymplectic quasi-Frobenius Lie superalgebra and let $\fl$ be an abelian Lie superalgebra. Let
$\fd:=\fl^*\oplus\fa\oplus\fl$. We define a bilinear form $\omega_\fd$ on $\fd$ by
\begin{equation}\label{orthoform}
    \omega_\fd(Z_1+a+L_1,Z_2+b+L_2):=Z_1(L_2)+\omega_\fa(a,b)-(-1)^{|Z_2||L_1|}Z_2(L_1),
\end{equation}
for all $a,b\in \fa$, for all $L_1,L_2\in \fl$ and all $Z_1,Z_2\in\fl^*$. Let us suppose that 
\begin{align}
    \xi(L)&\in\Der(\fa)~\forall L\in \fl;\label{qzz0}\\\label{qzz1}
    \ad\circ\alpha&=d_{\mathrm{CE}}^1\xi+\frac{1}{2}[\xi\wedge\xi]_\fa;\\\label{qzz3}
    d_{\xi}\alpha&=0;\\
    \gamma(L)\bigl([a,b]_\fa\bigl)&=\beta\bigl(\xi(L)(a),b\bigl)+(-1)^{|a||L|}\beta\bigl(a,\xi(L)(b)\bigl);\\\label{qzz5}
    \Ev(\gamma\wedge\alpha)&=0;\\\label{qzz2}
    d_{\widehat{\xi}}\gamma&=\widetilde{\beta}\circ \alpha;\\\label{qzz4}
    (-1)^{|L_1||L_3|}\epsilon(L_1,L_2)(L_3)&+(-1)^{|L_2||L_3|}\epsilon(L_3,L_1)(L_2)+(-1)^{|L_1||L_2|}\epsilon(L_2,L_3)(L_1)=0.   
\end{align}
Then the triple $(\fd,[\cdot,\cdot]_\fd,\omega_\fd)$ is an orthosymplectic quasi-Frobenius Lie superalgebra with the bracket $[\cdot,\cdot]_\fd$ defined as in \emph{Eq.} \eqref{extbracket} if and only if Eqs. \eqref{qzz0} through \eqref{qzz4} hold.
\end{Theorem}

\begin{proof}
The proof follows from the Jacobi identity for the bracket \eqref{extbracket} and the 2-cocycle condition for the bilinear form $\w_{\fd}$ defined in Eq. \eqref{orthoform}. Let us just prove Eqs. \eqref{qzz1} and \eqref{qzz2}. Let $a\in\fa$ and $L_1,L_2\in \fl$. The Jacobi identity gives 
\begin{equation}
    (-1)^{|a||L_2|}[a,[L_1,L_2]_\fd]_\fd+(-1)^{|a||L_1|}[L_1,[L_2,a]_\fd]_\fd+(-1)^{|L_2||L_1|}[L_2,[a,L_1]_\fd]_\fd=0.
\end{equation}
Therefore,
\begin{align*}
(-1)^{|a||L_2|}\Bigl( [a,\alpha(L_1,L_2)]_\fa+\beta(a,\alpha(L_1&,L_2)\Bigl)
+(-1)^{|a||L_1|}\Bigl(\gamma(L_1)(\xi(L_2)(a))+\xi(L_1)(\xi(L_2)(a))  \Bigl)\\
&-(-1)^{|L_1|(|a|+|L_2|)}\Bigl(\gamma(L_2)(\xi(L_1)(a))+\xi(L_2)(\xi(L_1)(a))\Bigl)\\
&=0.
\end{align*}
Upon collecting the coefficients in $\fl^*$, one obtains Eq. \eqref{qzz2}, since
$$d_{\widehat{\xi}}\gamma(L_1,L_2)(a)=-(-1)^{|L_1||L_2|}\gamma(L_2)\bigl(\xi(L_1)(a)\bigl)+\gamma(L_1)\bigl(\xi(L_2)(a) \bigl).$$
Upon collecting the coefficients in $\fa$, on obtains Eq. \eqref{qzz1}. Similarly, the other equations are obtained by carefully handling the $(-1)$ signs.
\end{proof}
\noindent\textbf{Remark.} In \cite{Fi}, pairs $(\xi,\alpha)$ satisfying Eqs \eqref{qzz0},\eqref{qzz1} and \eqref{qzz3} are called \textit{factor system}. See also \cite{AMR,Ho,N}, where factors systems are used to describe extensions of non-abelian Lie (super)algebras.

\subsection{A standard model for double extensions -- the periplectic case }\label{subsectionperi}   

Let $(\fa,[\cdot,\cdot]_\fa,\omega_\fa)$ be a  periplectic quasi-Frobenius Lie superalgebra and let $\fl$ be an abelian Lie superalgebra. We define a bracket $[\cdot,\cdot]_\fd$ on $\fd:=\Pi(\fl^*)\oplus\fa\oplus\fl$ by 
\begin{equation}\label{extbracketp}
\begin{array}{lclrcll}
[\Pi(\fl^*),\fd ] & = & 0; &   & [a,b]_\fd&=&\pi(\beta(a,b))+[a,b]_\fa; \\[2mm] 
 [L,a]_\fd & = & \pi(\gamma(L)(a))+\xi(L)(a);& & [L_1,L_2]_\fd&=&\pi(\epsilon(L_1,L_2))+\alpha(L_1,L_2),
 \end{array}
\end{equation}
for all $a,b\in\fa$ and for all $L,L_1,L_2\in\fl$, where the maps $\alpha,\beta,\gamma,\xi,\epsilon$ are defined as in Eqs. \eqref{beta} to \eqref{alpha}. We define a bilinear form on $\fd$ by
\begin{equation}\label{perform}
    \omega_\fd\bigl(\pi(Z_1)+a+L_1,\pi(Z_2)+b+L_2\bigl):=Z_1(L_2)+\omega_\fa(a,b)-(-1)^{(|Z_2|+1)|L_1|}Z_2(L_1),
\end{equation}
for all $a,b\in \fa$, for all $L_1,L_2\in \fl$ and for all $\pi(Z_1),\pi(Z_2)\in\Pi(\fl^*)$.

\sssbegin{Theorem}[Superization of \text{\cite[Theorem 2]{Fi}}]\label{peristandard}
The triple $(\fd,[\cdot,\cdot]_\fd,\omega_\fd)$ is a periplectic quasi-Frobenius Lie superalgebra with the bracket $[\cdot,\cdot]_\fd$ defined in \emph{Eq.} \eqref{extbracketp} if and only if Eqs. \eqref{qzz0} through  \eqref{qzz4} hold.
\end{Theorem}

\begin{proof}
The proof follows as in Thm. \ref{orthostandard}.
\end{proof}

\subsection{Equivalence to the standard model}    

Let $(\fa,[\cdot,\cdot]_\fa,\w_{\fa})$ be a quasi-Frobenius Lie superalgebra and let $\fl$ be an abelian Lie superalgebra. Let $(\fg,\fj,i,p)$ be a double extension of $\fa$ by $\fl$ (see Def. \ref{doublextdef}). Since $\fj$ is isotropic, there exists an isotropic subspace $V_\fl\subset\fg$ such that $\fg=\fj^\perp\oplus V_\fl$.

\sssbegin{Lemma}
Let $V_\fa:=(\fj\oplus V_\fl)^\perp$. We have the decomposition $\fg=V_{\fa}\oplus\fj\oplus V_{\fl}$.
\end{Lemma}

\begin{proof}
   Let $j+v\in(\fj\oplus V_\fl)\cap(\fj^\perp\cap V_\fl^\perp)$. Since $j+v\in\fj^\perp$, and $\fj^\perp\cap V_\fl=0$, it follows that $v=0$. Similarly, we obtain $j=0$. Therefore, $V_\fa\cap(\fj\oplus V_\fl)=(\fj\oplus V_\fl)\cap(\fj^\perp\cap V_\fl^\perp)=0$. As a result, we obtain the desired decomposition, since $\fg=V_\fa+(\fj\oplus V_\fl)$. 
\end{proof}
Here is the main Theorem of this paper. 
\sssbegin{Theorem}[Superization of \text{\cite[Theorem 4]{Fi}}]\label{equivstandard} Let $\fl$ be an abelian Lie superalgebra, and let $(\fa,[\cdot,\cdot]_\fa,\w_{\fa})$ be a  quasi-Frobenius Lie superalgebra such that $\omega_\fg$ is homogeneous. Let $(\fg,\fj,i,p)$ be a double extension of $\fa$ by $\fl$.
\begin{enumerate}
   \item \underline{The orthosymplectic case: $|\omega_\fa|=\ev$}.
The symplectic double extension  $(\fg,\fj,i,p)$ is equivalent to the standard model $(\fd,[\cdot,\cdot]_\fd,\omega_\fd)$ as described in Theorem \ref{orthostandard}.

    \item \underline{The periplectic case $|\omega_\fa|=\od$}.
The symplectic double extension  $(\fg,\fj,i,p)$ is equivalent to the standard model $(\fd,[\cdot,\cdot]_\fd,\omega_\fd)$ as described in Theorem \ref{peristandard}.
\end{enumerate}
\end{Theorem}

\begin{proof}
We will only prove Thm. \ref{equivstandard} in the periplectic case. Let $(\fa,[\cdot,\cdot]_\fa,\w_{\fa})$ be a periplectic quasi-Frobenius Lie superalgebra. Let $(\fg,\fj,i,p)$ be a double extension of $\fa$ by $\fl$. Recall that the periplectic form on $\fg$ is given by  Eq. \eqref{perform}.  Let $\sigma:\fg\rightarrow \fg/\fj$ be the natural projection and let $\tilde{p}:\fg\rightarrow\fl$ be the composition $\tilde{p}=p\circ\sigma$. Let $s: \fl\rightarrow V_{\fl}$ be a linear isomorphism such that $\tilde{p}\circ s=\id$. Recall that the map $i:\fa\rightarrow\fj^{\perp}/\fj$ is an isomorphism of quasi-Frobenius Lie superalgebras. There is a map $t:\fa\rightarrow V_{\fa}$ defined by $$i(a)=\sigma(t(a))=t(a)+\fj\in \fj^{\perp}/\fj.$$
Moreover, we define a map $p^*:\Pi(\fl^*)\rightarrow\fg$ by
$$\w_{\fg}\bigl(p^*(\pi(Z)),s(L)\bigl)=Z(\tilde{p}\circ s)(L)=Z(L),~\forall Z\in \fl^*,~\forall L\in\fl.    $$
Since $\w_{\fg}$ is non-degenerate, the map $p^*$ defines an isomorphism $\Pi(\fl^*)\rightarrow\fj$. The following diagram sums up the situation.

\begin{equation}
\begin{tikzcd}
            &                                                   & 0 \arrow[d]                                               &                                                      &   \\
            &                                                   & V_{\mathfrak{l}} \arrow[d, "\subset" description]         &                                                      &   \\
0 \arrow[r] & V_{\mathfrak{a}} \arrow[r, "\subset" description] & \mathfrak{g} \arrow[d, "\sigma"'] \arrow[rd, "\tilde{p}"] &                                                      &   \\
0 \arrow[r] & \mathfrak{a} \arrow[r, "i"'] \arrow[u, "t"]       & \mathfrak{g}/\mathfrak{j} \arrow[r, "p"']                 & \mathfrak{l} \arrow[r] \arrow[luu, "s"', bend right] & 0
\end{tikzcd}
\end{equation}
For all $a,b\in\fa$, we define
\begin{equation}
    \w_\fa(a,b)=\w_{\fj^{\perp}/\fj}(i(a),i(b))=\w_\fg(t(a),t(b)).
\end{equation}
The proof then consists of four steps. (Recall that $\fd=\Pi(\fl^*)\oplus\fa\oplus\fl$.)\\

\noindent\textbf{\underline{Step $1$.}} We will be defining a bracket  $[\cdot,\cdot]_{\fd}:\fd\times\fd\rightarrow\fd$ and a symplectic form $\w_\fd:\fd\times\fd\rightarrow\K$. Let
$$\xi:~\fl\rightarrow\End(\fa);~ \gamma:~\fl\rightarrow\Hom(\fa,\fl^*);~\text{and } \epsilon:~\fl\times\fl\rightarrow\fl^*$$
defined for all $a\in\fa$ and for all $L,L_1,L_2\in\fl$ by
\begin{align}
    \label{onze}i(\xi(L)(a))&:=[s(L),t(a)]_\fg+\fj\in \fj^{\perp}/\fj;\\\label{douze}
    p^*(\gamma(L))(a)&:=[s(L),t(a)]_\fg-t\bigl(\xi(L)(a)\bigl);\\\label{treize}
    \epsilon(L_1,L_2)&:=\w_{\fg}\bigl([s(L_1),s(L_2)]_\fg,s(\cdot)\bigl).
\end{align}
Since $\w_\fg$ is closed, it follows that $[s(L),t(a)]_\fg\in\fj^{\perp}$. Therefore, $\xi$ is a well-defined map. Moreover, $\sigma\bigl([s(L),t(a)]_{\fg}\bigl)=[s(L),t(a)]_{\fg}+\fj$ and 
\[
\sigma\circ t(\xi(L)(a))=i(\xi(L)(a))=\sigma\bigl([s(L),t(a)]_{\fg}).
\] Thus, we have $[s(L),t(a)]_\fg-t\bigl(\xi(L)(a)\bigl)\in\fj$ and that $\gamma$ is well-defined. The maps $\beta$ and $\alpha$ can be defined as in Eqs. \eqref{beta} and \eqref{alpha}, by 
\begin{align}
        \beta(a,b)(L)&=-(-1)^{|L||b|}\w_\fa\bigl(a,\xi(L)(b)\bigl)+(-1)^{|a|(|L|+|b|)}\w_\fa\bigl(b,\xi(L)(a)\bigl)\\
        \w_{\fa}(a,\alpha(L_1,L_2))&=(-1)^{|L_2|(|a|+|L_1|)}\gamma(L_2)(a)(L_1)-(-1)^{|a||L_1|}\gamma(L_1)(a)(L_2),
\end{align}
for all $a,b\in \fa$ and all $L_1,L_2,L\in \fl$.
Now, the bracket $[\cdot,\cdot]_\fd$ and the symplectic form $\w_\fd$ can be defined as in Eqs. \eqref{extbracketp} and \eqref{perform}.\\

\noindent\textbf{\underline{Step $2$.}} We define a linear isomorphism of vector spaces
\begin{align}
    \phi:\fd=\Pi(\fl^*)\oplus\fa\oplus\fl&\rightarrow \fg\\\nonumber
                (\pi(Z)+a+L)&\mapsto\bigl( p^*(\pi(Z))+t(a)+s(L) \bigl).
\end{align}

\noindent\textbf{\underline{Step $3$.}} We shall show that $\phi\circ[\cdot,\cdot]_\fd=[\phi(\cdot),\phi(\cdot)]_\fg$ and that $\w_\fd=\phi^*\w_\fg$.
We need the two following lemmas.
\sssbegin{Lemma}\label{ab}
    Let $a,b\in\fa$. Then $[t(a),t(b)]_\fg-t([a,b]_\fa)\in\fj$. Moreover, 
    \[p^*\bigl(\pi(\beta(a,b))\bigl)=[t(a),t(b)]-t([a,b]_\fa),\quad  \text{ and } \quad \phi\circ[a,b]_\fd=[\phi(a),\phi(b)]_\fg.
    \]
\end{Lemma}

\begin{proof}
Let $a,b\in\fa$ and $L\in\fl$. 
   $$\sigma\bigl( [t(a),t(b)]_\fg-t([a,b]_\fa) \bigl)=[i(a),i(b)]_{\fj^{\perp}/\fj}-i\bigl([a,b]_\fa\bigl)=0.$$ Therefore, there exists $\pi(X)\in\Pi(\fl^*)$ such that $[t(a),t(b)]_\fg-t([a,b]_\fa)=p^*(\pi(X))$. A careful computation shows that
   $$p^*\bigl(\pi(\beta(a,b))  \bigl)=p^*\Bigl(\pi\bigl(\w_\fg([t(a),t(b)],s(\cdot))\bigl)\Bigl).$$ Therefore, 
   \begin{align*}
\w_\fg\Bigl(p^*\bigl(\pi(\beta(a,b))\bigl),s(L)\Bigl)&=\w_\fg\Bigl(p^*\Bigl(\pi\bigl(\w_\fg([t(a),t(b)],s(L))\bigl)\Bigl)\Bigl).\\
&=\w_\fg\bigl([t(a),t(b)],s(L)\bigl)\\
&=\w_\fg\bigl([t(a),t(b)]-t([a,b]_\fa),s(L)\bigl),
\end{align*}
since $t([a,b]_\fa)$ and $s(L)$ are orthogonal. Since $\w_\fg(\fj,\fj\oplus V_\fa)=0$ and since $s$ is an isomorphism, we have  $$\w_\fg\Bigl(p^*\bigl(\pi(\beta(a,b))\bigl),\fg\Bigl)=\w_\fg\bigl([t(a),t(b)]-t([a,b]_\fa),\fg\bigl).$$ Since $\omega_\fg$ is non-degenerate, it follows that   $p^*\bigl(\pi(\beta(a,b))\bigl)=[t(a),t(b)]-t([a,b]_\fa)$. Thus,  $\pi(X)=\pi(\beta(a,b))$.  Therefore,
$$\phi([a,b]_\fd)=p^*\Bigl(\pi\bigl(\beta(a,b)\bigl)\Bigl)+t\bigl([a,b]_\fa\bigl)=[t(a),t(b)]_\fg=[\phi(a),\phi(b)]_\fg,~\forall a,b\in\fa.$$\end{proof}

\sssbegin{Lemma}\label{LL}
        Let $L_1,L_2\in\fl$. Then,
        $$t\bigl(\alpha(L_1,L_2)\bigl)=\bigl[s(L_1),s(L_2)\bigl]_\fg-p^*\Bigl(\pi\bigl(\epsilon(L_1,L_2)\bigl)\Bigl).   $$
\end{Lemma}
\begin{proof}
    Let $a\in\fa$ and $L_1,L_2\in\fl$.\\
    
\textbf{Claim.} We have $[s(L_1),s(L_2)]_\fg-p^*\bigl(\pi(\epsilon(L_1,L_2))\bigl)\in V_\fa=\fj^{\perp}\cap V_\fl^{\perp}$. 

\noindent\textbf{Proof of the claim.} First, we have
\begin{align*}
    \w_g\Bigl([s(L_1),s(L_2)]_\fg-p^*\bigl(\pi(\epsilon(L_1,L_2))\bigl),s(L)  \Bigl)
    &= \w_\fg\bigl([s(L_1),s(L_2)]_\fg,s(L)\bigl)-\epsilon(L_1,L_2)(L)\\&=0~~ (\text{Eq}.~\eqref{treize});
\end{align*}
Therefore, $[s(L_1),s(L_2)]_\fg-p^*\bigl(\pi(\epsilon(L_1,L_2))\bigl)\in V_\fl^{\perp}$. Moreover,
$$\w_g\Bigl([s(L_1),s(L_2)]_\fg-p^*\bigl(\pi(\epsilon(L_1,L_2))\bigl),\fj \Bigl)=0 $$ using the $2$-cocycle condition. Hence the claim. Then,
    \begin{align*}       \w_\fg\bigl(t(a),t(\alpha(L_1,L_2))\bigl)&=(-1)^{|L_2|(|a|+|L_1|)}\gamma(L_2)(a)(L_1)-(-1)^{|a||L_1|}\gamma(L_1)(a)(L_2)\\
        &=(-1)^{|L_2|(|a|+|L_1|)}\w_\fg\Bigl(p^*\circ\pi(\gamma(L_2)(a)),s(L_1)   \Bigl)\\&~~~-(-1)^{|a||L_1|}\w_\fg\Bigl(p^*\circ\pi(\gamma(L_1)(a)),s(L_2)   \Bigl)\\
        (\text{using Eq. }\eqref{douze})~~~~~&=(-1)^{|L_2|(|a|+|L_1|)}\w_{\fg}\Bigl([s(L_2),t(a)]_\fg-t\bigl(\xi(L_2)(a)\bigl),s(L_1)\Bigl)\\
        &~~~-(-1)^{|a||L_1|}\w_\fg\Bigl( [s(L_1),t(a)]_\fg-t\bigl(\xi(L_1)(a)\bigl),s(L_2)\Bigl)\\
        &=(-1)^{|L_2|(|a|+|L_1|)}\w_\fg\bigl([s(L_2),t(a)]_\fg,s(L_1)\bigl)+ \w_\fg\bigl([s(L_1),t(a)]_\fg,s(L_2)\bigl)\\   
        &=\w_\fg\bigl(t(a),[s(L_1),s(L_2)]_\fg  \bigl).
    \end{align*}

We deduce that 
 \begin{equation}     \w_\fg\Bigl(t(a),t\bigl(\alpha(L_1,L_2)\bigl)\Bigl)=\w_\fg\Bigl(t(a), [s(L_1),s(L_2)]_\fg-p^*\bigl(\pi(\epsilon(L_1,L_2))\bigl)   \Bigl),
 \end{equation}
 since $p^*\circ\pi(\epsilon(L_1,L_2))$ is orthogonal to $V_\fa$. Using the claim and using the fact that $\w_\fg$ is non-degenerate, we finally have
  \begin{equation}     t\bigl(\alpha(L_1,L_2)\bigl)= [s(L_1),s(L_2)]_\fg-p^*\bigl(\pi(\epsilon(L_1,L_2))\bigl),
 \end{equation}\end{proof}
Let us finish the proof of Step 3. Let $a\in\fa$ and $L\in\fl$. Using Eq. \eqref{douze}, we have 
\begin{equation}\label{La}
    \phi([L,a]_\fd)=p^*\bigl(\pi(\gamma(L)(a))\bigl)+t(\xi(L)(a))=[s(L),t(a)]_\fg=[\phi(L),\phi(a)]_\fg.
\end{equation}

Lemmas \ref{ab} and \ref{LL} together with Eq. \eqref{La} implies that $\phi\circ[\cdot,\cdot]_\fd=[\phi(\cdot),\phi(\cdot)]_\fg$. To finish Step 3, it remains to show that $\w_\fd=\phi^*\w_\fg$. This is an immediate consequence of the definition of the form $\w_\fd$, see Eq. \eqref{perform}.\\

\noindent\textbf{\underline{Step $4$.}} Step 3 ensures that $(\fd,[\cdot,\cdot]_\fd,\omega_\fd)$ is a periplectic Lie superalgebra, isomorphic to $(\fg,[\cdot,\cdot]_\fg,\w_\fg)$. We show that it is equivalent to $(\fg,[\cdot,\cdot]_\fg,\w_\fg)$ in the sense of Definition \ref{eqextt}. Consider the diagram
\begin{equation}\label{encoreundiagram}
\begin{tikzcd}
            &                                          & \mathfrak{d}/\Pi(\mathfrak{l}^*) \arrow[rd] \arrow[dd, "\bar{\phi}"] &                        &   \\
0 \arrow[r] & \mathfrak{a} \arrow[rd, "i"'] \arrow[ru] &                                                                      & \mathfrak{l} \arrow[r] & 0. \\
            &                                          & \mathfrak{g}/\mathfrak{j} \arrow[ru, "p"']                           &                        &  
\end{tikzcd}
\end{equation}
Let $a\in\fa$ and $L\in\fl$. Let $\Tilde{p}=p\circ\sigma$, where $\sigma$ is the natural projection $\fg\rightarrow \fg/\fj$. Since $\Bar{\phi}(a)=t(a)+\fj$, we have $i(a)=\sigma(t(a))=\Bar{\phi}(a)$. Moreover, \[
p(\bar{\phi}(L))=p(\sigma\circ s(L))=\Tilde{p}(s(L))=L.\] Therefore, the diagram \eqref{encoreundiagram} commutes and the map $\phi$ defines an equivalence of double extensions.
\end{proof}

\subsection{Equivalence classes} In this section, we describe a general procedure to compute equivalence classes of doubles extensions. Let $(\fg,\w_\fg)$ be a symplectic Lie superalgebra and let $\fl$ be an abelian superalgebra. Given a linear map $\tau:\fl\rightarrow\fa$, we define
\begin{equation}
    \tau^*:\fa\rightarrow\fl^*,~~a\mapsto-\w_\fa\bigl(a,\tau(\cdot)\bigl).
\end{equation}
    
\sssbegin{Proposition}[Superization of \text{\cite[Lemma 15]{Fi}}]
    Let $\fd_1$ \textup{(}resp. $\fd_2$\textup{)} be a symplectic  double extension of  $(\fg,\w_\fg)$ by $\fl$ given by the maps $\xi_1, \gamma_1, \epsilon_1$ \textup{(}resp. $\xi_2, \gamma_2, \epsilon_2$\textup{)}  \textup{(}see Thm. \ref{orthostandard}\textup{)}. Then, the double extensions $\fd_1$ and $\fd_2$ are equivalent if and only if there exists a linear map $\tau:\fl\rightarrow\fa$ such that
    \begin{align}
        \xi_2&=\xi_1-\ad\circ\tau:\\\label{bbbb}
        \gamma_2&=\gamma_1+\tau^*\circ\xi_1-\Tilde{\beta_1}\circ\tau-\tau^*\circ\ad\circ\tau;\\\label{Emm}
        \epsilon_2&=\epsilon_1+\tau^*\circ\alpha_1-\Ev(\gamma_1\wedge\tau)+\frac{1}{2}\beta_1(\tau\wedge\tau)+\tau^*d^1_{\xi_1}\tau, 
    \end{align}
    where $\beta_1$ and $\alpha_1$ are given as in Eqs. \eqref{beta},\eqref{alpha} using the maps $\gamma_1$ and $\epsilon_1$.
\end{Proposition}
\begin{proof}
Our proof will be limited to the periplectic case. Suppose that $(\fd_1,\w_1)$ and $(\fd_2,\w_2)$ are equivalent doubles extensions of $(\fa,[\cdot,\cdot ]_\fa, \w_\fa)$ by $\fl$ by an isomorphism $\phi:\fd_1\rightarrow\fd_2$ (see Definition \ref{eqextt}). Since $\overline{\phi}\circ i_1=i_2$ and $p_2\circ\overline{\phi}=p_1$, it follows that $\phi(a)\in\fa\oplus\Pi(\fl^*)$ and that $\phi|_\fl=\id$. Therefore, $\phi$ is given with respect to the decomposition $\fg=\Pi(\fl^*)\oplus\fa\oplus\fl$ by a matrix of the form
\begin{equation}
    \phi=\begin{pmatrix}\id&c&b\\0&\id&\tau\\0&0&\id\end{pmatrix},
\end{equation}
where $b:\fl\rightarrow \Pi(\fl^*)$, $c:\fa\rightarrow \Pi(\fl^*)$ and $\psi:\Pi(\fl^*)\rightarrow \Pi(\fl^*)$ are linear maps. From $\varphi^*\w_1=\w_2$, we obtain $c=\pi\circ\tau^*$ and
\begin{equation}
        \tau^*\circ\tau(L)(L')=\pi\bigl(b(L)\bigl)(L')-(-1)^{|L||L'|}\pi\bigl(b(L')\bigl)(L),~~\forall L,L'\in\fl.
\end{equation}
Let $a_1,a_2\in\fa$. Since $\phi([a_1,a_2]_1)=[\phi(a_1),\phi(a_2)]_2$, we have
\begin{equation}\label{bbb}
    \beta_2(a_1,a_2)=\beta_1(a_1,a_2)+\tau^*([a_1,a_2]_\fa).
\end{equation}
Moreover, expanding $\phi([L,a])=[\phi(L),\phi(a)]_2$ for all $a\in\fa,~L\in\fl$ gives
\begin{align}
    \xi_2(L)(a)&=\xi_2(L)(a)-[\tau(L),a]_\fa\\\label{bb}
    \beta_2\bigl(\tau(L),a\bigl)+\gamma_2(L)(a)&=\gamma_1(L)(a)+\tau^*\circ\xi_1(L)(a).
\end{align}
Eqs. \eqref{bbb} and \eqref{bb} together give \eqref{bbbb}. Let $L_1,L_2\in\fl$. We have
\begin{align}
    [\phi(L_1),\phi(L_2)]_2&=[\tau(L_1),\tau(L_2)]+\pi\bigl(\beta_2(\tau(L_1),\tau(L_2))\bigl)\\\nonumber
    &-(-1)^{|L_1||L_2|}\Bigl(\pi\bigl(\gamma_2(L_2)(\tau(L_1))\bigl)+\xi_2(L_2)(\tau(L_1))\Bigl)\\\nonumber
    &+\pi\bigl(\gamma_2(\tau(L_2))\bigl)+\xi_2(L_1)\bigl(\tau(L_2)\bigl)\\\nonumber
    &+\pi\bigl(\epsilon_2(L_1,L_2)\bigl)+\alpha_2(L_1,L_2).
\end{align}
Moreover, 
\begin{equation}
    \phi\bigl([L_1,L_2]_1\bigl)=\pi\bigl(\epsilon_1(L_1,L_2)\bigl)+\pi\bigl(\tau^*(\alpha_1(L_1,L_2))\bigl)+\alpha_1(L_1,L_2).
\end{equation}
Using \eqref{bbbb} and since $\phi\bigl([L_1,L_2]_1\bigl)=[\phi(L_1),\phi(L_2)]_2$, one gets
\begin{align}\label{Em}
    \beta_2\bigl(\tau(L_1),\tau(L_2)\bigl)-&(-1)^{|L_1||L_2|}\gamma_2(L_2)(\tau(L_1))+\gamma_2(L_1)(\tau(L_2))+\epsilon_2(L_1,L_2)\\\nonumber&=\epsilon_1(L_1,L_2)+\tau^*(\alpha_1(L_1,L_2)).
\end{align}
Using \eqref{bbb} and \eqref{bb}, Eq. \eqref{Em} is equivalent to \eqref{Emm}.\end{proof}

\section{Examples of symplectic doubles extensions}\label{sectionexamples}

The purpose of this section is to compute explicitly examples of doubles extensions in the orthosymplectic and periplectic cases, respectively. In our examples, we use at least a two-dimensional ideal, and the Lie superalgebras are taken from \cite{Ba}. For various examples where the ideal is one-dimensional, see \cite{BM, BEM}.
\subsection{The Lie superalgebra $C^3+A$, see \cite{Ba}} Consider the Lie superalgebra $\fa=C^3+A$ with the brackets given in the basis (even $|$ odd) $e_1,e_2|e_3,e_4$ by \begin{equation}[e_1,e_4]=e_3;~~[e_4,e_4]=e_2.\end{equation}
As shown in \cite{BM}, $\fa$ is orthosymplectic quasi-Frobenius with the form given by\footnote{As in \cite{BM,BEM}, we adopt the following convention: $\langle e_i^*, e_j\rangle=\delta_{ij}$ and  $\langle e_i^*\otimes e_j^*, e_k\otimes e_l\rangle= (-1)^{|e_k||e_j|} \langle e_i^*, e_k \rangle\langle e_j^*, e_l\rangle $.} 
\begin{equation}
\omega=2e_1^*\wedge e_2^* - e_3^*\wedge e_4^*.
\end{equation}
Let $\fl$ be an abelian Lie superalgebra of superdimension $(1|1)$ with basis $L_1|L_2$. We consider the two derivations of $C^3+A$ given by
\begin{equation}
d_1:=e_2\otimes e_1^*;~~d_2:=e_2\otimes e_4^*.
\end{equation}
We define a map $\xi:\fl\rightarrow\Der(\fa)$ by $\xi(L_1)=d_1$ and $\xi(L_2)=d_2$.
Using Eq. \eqref{beta}, we obtain
\begin{equation}
    \beta=2L_2^*\otimes(e_1^*\wedge e_4^*).
\end{equation}
Moreover, consider
\begin{equation}
    \gamma=L_1^*\otimes e_4^*\otimes L_2^*-L_2^*\otimes e_4^*\otimes L_1^*+L_1^*\otimes e_1^*\otimes L_2^*.
\end{equation}
Let  $\alpha:\fl\times\fl\rightarrow\fa$ be a bilinear map defined as in Eq. \eqref{alpha}. Then, Eqs. \eqref{qzz3},\eqref{qzz5} and \eqref{qzz2}  are trivially satisfied and Eq. \eqref{qzz1} gives $\alpha(L,L')\in \Span(e_2|e_3)$, $\forall L,L'\in\fl$. Furthermore, Eq. \eqref{alpha} gives
\begin{equation}
    \alpha=\frac{1}{2}e_2\otimes(L_2^*\wedge L_2^*).
\end{equation}
Finally, let
\begin{equation}
    \epsilon:=L_2^*\wedge(L_1^*\wedge L_2^*).
\end{equation}
A direct comutation shows that the the map $\epsilon$ satisfies Eq. \eqref{qzz4}.  

The Lie superalgebra structure on $\fl^*\oplus \fa\oplus \fl$ is then given by

\begin{equation}
\begin{array}{lllllllll}
[e_1,e_4 ] & = & e_3+2L_2^*; & [e_1,L_1]&=&-e_2;  & [e_1,L_2]&=&-2L_2^*;   \\[2mm] 
[e_4,e_4 ] & = & e_2; & [e_4,L_1]&=&L^*_2;  & [e_4,L_2]&=&e_2+L_1^*;   \\[2mm] 
[L_1,L_2 ] & = & L_2^*; & [L_2,L_2]&=&-e_2.  & &&   \\[2mm] 
 \end{array}
\end{equation}
It is worth noting that this Lie  superalgebra is $2$-step nilpotent.

\subsection{The Lie superalgebra $C^1_{1/2}+A$, see \cite{Ba}} Consider the Lie superalgebra $\fa=C^1_{1/2}+~A$ with the brackets given in the basis (even $|$ odd) $e_1,e_2|e_3,e_4$ by 
\begin{equation}
[e_1,e_2]=e_2;~~[e_1,e_3]=\frac{1}{2}e_3;~~[e_3,e_3]=e_2.\end{equation}
As shown in \cite{BM}, $\fa$ is orthosymplectic quasi-Frobenius with the form given by 
\begin{equation}
\omega=e_1^*\wedge e_2^*-\frac{1}{2}e_3^*\wedge e_3^*-\frac{1}{2}e_4^*\wedge e_4^*.
\end{equation}
Let $\fl$ be the abelian Lie superalgebra of superdimension $(1|1)$ with basis $L_1|L_2$. We consider the two derivations of $C^1_{1/2}+A$ given by
\begin{equation}
d_1:=\ad_{e_1}+\ad_{e_2};~~d_2:=e_4\otimes e_1^*.
\end{equation}
We define a map $\xi:\fl\rightarrow\Der(\fa)$ by $\xi(L_1)=d_1$ and $\xi(L_2)=d_2$.
Using Eq. \eqref{beta}, we obtain
\begin{equation}
    \beta=-L_1^*\otimes\Bigl(e_1^*\wedge e_2^*+\frac{1}{2}e_3^*\wedge e_3^*\Bigl)+L_2^*\otimes (e_1^*\wedge e_4^*).
\end{equation}
Moreover, consider
\begin{equation}
    \gamma=L_1^*\otimes e_2^*\otimes L_1^*, \quad \alpha\equiv 0, \quad \epsilon=L_1^*\wedge(L_2^*\wedge L_1^*).
\end{equation}
The map $\epsilon$ satisfies Eq. \eqref{qzz4}. 

The Lie superalgebra structure on $\fl^*\oplus \fa\oplus \fl$ is then given by

\begin{equation}
\begin{array}{lllllllll}
[e_1,e_2 ] & = & e_2-L_1^*; & [e_1,e_3]&=&\frac{1}{2}e_3;  & [e_1,e_4]&=&L_2^*;   \\[2mm] 
[e_1,L_1 ] & = & e_2; & [e_1,L_2]&=&-e_2;  & [e_2,L_1]&=&-e_2+L_1^*;   \\[2mm] 
[e_3,e_3 ] & = & e_2-L_1^*; & [e_3,L_1]&=&-\frac{1}{2}e_3;  & [L_1,L_2]&=&-L_1^*.   \\[2mm] 
 \end{array}
\end{equation}

\subsection{The Lie superalgebra $(2A_{1,1}+2A)^3_{1/2}$, see \cite{Ba}} Consider the Lie superalgebra $\fa=(2A_{1,1}+2A)^3_{1/2}$ with the brackets given in the basis (even $|$ odd) $e_1,e_2|e_3,e_4$ by \begin{equation}[e_3,e_3]=e_1;~~[e_3,e_4]=\frac{1}{2}(e_1+e_2);~~[e_4,e_4]=e_2.\end{equation}
It has been shown in \cite{BM} that $\fa$ is periplectic quasi-Frobenius with the form given by 
\begin{equation}
\omega=e_2^*\wedge e_3^* - e_1^*\wedge e_4^*.
\end{equation}
Let $\fl$ be the abelian Lie superalgebra of superdimension $(0|2)$ with basis $0|L_1,L_2$. We consider the two odd derivations of $(2A_{1,1}+2A)^3_{1/2}$ given by
\begin{equation}
d_1:=e_1\otimes e_3^*;~~d_2:=e_1\otimes e_4^*.
\end{equation}
We define a map $\xi:\fl\rightarrow\Der(\fa)$ by $\xi(L_1)=d_1$ and $\xi(L_2)=d_2$.
Using Eq. \eqref{beta}, we obtain
\begin{equation}
    \beta=L_1^*\otimes(e_3^*\wedge e_4^*)+L_2^*\otimes(e_4^*\wedge e_4^*).
\end{equation}
Moreover, consider
\begin{equation}
    \gamma=L_1^*\otimes e_3^*\otimes L_2^*+L_2^*\otimes e_3^*\otimes L_1^*.
\end{equation}
Let $\alpha:\fl\times\fl\rightarrow\fa$ be a bilinear map defined as in Eq. \eqref{alpha}. Then, Eqs. \eqref{qzz3},\eqref{qzz5} and \eqref{qzz2} are trivially satisfied and Eq. \eqref{qzz1} gives $\alpha(L,L')\in \Span(e_1,e_2|0)$, $\forall L,L'\in\fl$. Furthermore, Eq. \eqref{alpha} gives
\begin{equation}
    \alpha=-2e_2\otimes(L_1^*\wedge L_2^*).
\end{equation}
Finally, let
\begin{equation}
    \epsilon=L_2^*\wedge(L_1^*\wedge L_2^*)-L_1^*\wedge(L_2^*\wedge L_2^*).
\end{equation}
A direct computation shows that  $\epsilon$ satisfies Eq. \eqref{qzz4}.  

The Lie superalgebra structure on $\Pi(\fl^*)\oplus \fa\oplus \fl$ is then given by

\begin{equation}
\begin{array}{lllllllll}
[e_3,e_3 ] & = & e_1; & [e_3,e_4]&=&\frac{1}{2}(e_1+e_2)-\pi(L_1^*);  &&&   \\[2mm] 
[e_4,e_4]&=&e_2-2\pi(L_2^*);&[L_1,e_3 ] & = & e_1+\pi(L_2^*); & [L_2,e_3]&=&\pi(L_1^*);      \\[2mm] 
[L_2,e_4]&=&e_1;&[L_1,L_2 ] & = & -2e_2+\pi(L_2^*); & [L_1,L_2]&=&-\pi(L_1^*). \\[2mm] 
 \end{array}
\end{equation}
\section{Open problems}
A list of open problems is given below.

\begin{enumerate}[(i)]

\item Describe symplectic double extensions in the case where the center is trivial. For a 1-dimensional ideal, this problem has already been solved in \cite{BM}.

\item Over a field of characteristic $p$, if the Lie superalgebra $\fa$ is restricted, under which conditions its double extension is restricted as well? The 1-dimensional case has already been studied in \cite{BEM}.

\item\label{iii} In \cite[Corollaire 21]{H}, the author proved that if $\fg$ is a left-symmetric finite-dimensional Lie algebra, then $[\fg,\fg]\neq\fg$ (see also \cite[Thms 6 and 7]{C}). Since quasi-Frobenius Lie algebras admit a compatible left-symmetric structure  (see e.g. \cite{DM2}), they do satisfy this condition. Is it also true for quasi-Frobenius Lie superalgebras?

\item The existence of an isotropic ideal $\fj$ is crucial to the process  of double extensions. In \cite[Lemme 1.1]{DM2}, Dardi\'e and Medina proved  that any quasi-Frobenius Lie algebra admits a non trivial abelian ideal $\fj$, which is either isotropic or symplectic. The proof relies on \cite[Corollaire 21]{H} that we discussed in point \eqref{iii}. In the case where the quasi-Frobenius Lie algebra does not admit any isotropic ideal, they showed that $\fg =[\fg,\fg] \rtimes [\fg,\fg]^{\perp}$, a semi direct product of the abelian derived ideal $[\fg,\fg]$ by $[\fg,\fg]^{\perp}$ (\cite[Th\'eor\`eme 1.3]{DM2}). It would be very interesting to investigate this result in the super settings.

\item Classify quasi-Frobenius Lie superalgebras of dimension 6 based on the same ideas as in \cite{Fi}.\\

\end{enumerate}

\noindent\textbf{Acknowledgements.} We would like to thank S. Benayadi for several stimulating and fruitful discussions.


\end{document}